\documentclass[a4paper,10pt]{article}
\usepackage[utf8]{inputenc}
\usepackage{amsmath,amsthm,amssymb}
\usepackage{amsfonts}
 
 \theoremstyle{plain}
 \newtheorem{lem}{Lemma}
 \newtheorem*{lem*}{Lemma}
\newtheorem{cor}{Corollary}
\newtheorem{thm}{Theorem}
 
\theoremstyle{definition}
\newtheorem{defn}{Definition}

\theoremstyle{definition} 
\newtheorem*{defn*}{Definition}
 
\theoremstyle{definition}
\newtheorem{rmk}{Remark}

\newtheorem{conj}{Conjecture}

\newtheoremstyle{named}{}{}{\itshape}{}{\bfseries}{.}{.5em}{\thmnote{#3}}
\theoremstyle{named}

\newtheoremstyle{namedproof}{}{}{}{}{\itshape}{.}{.5em}{\thmnote{#3}}
\theoremstyle{namedproof}
\newtheorem*{namedproof}{namedproof}

\newcommand{\bbb}[1]{\textbf{#1}}

\newcommand{\wt}[1]{\widetilde{#1}}

\newcommand{\mS}{\mathcal{S}}
\newcommand{\comm}[1]{}
\newcommand{\wh}[1]{\widehat{#1}}

\renewcommand{\d}{\mathrm{d}}
\renewcommand{\S}{\mathrm{Sym}}

\title{On hereditarily just infinite profinite groups obtained via iterated wreath products}
\author{MATTEO VANNACCI}
\begin{document}

\maketitle

\begin{abstract}
\noindent We study a generalisation of the family of non-(virtually pro-$p$) hereditarily just infinite profinite groups introduced by J.\! S.\! Wilson in 2010. We prove that this family contains groups of finite lower rank. We also show that many groups in this family are not topologically finitely presentable. 
\end{abstract}



\section{Introduction} 

Let $(m_k)_{k\in \mathbb{N}}$ be a sequence of positive integers. Throughout this paper we denote by $\mS$ a sequence of finite transitive permutation groups $(S_k \le \S(m_k))_{k\in \mathbb{N}}$ such that each $S_k$ is non-abelian simple as an abstract group. 
\begin{defn}
 We say that a profinite group $G$ is just infinite if it is infinite, and every non-trivial closed normal subgroup is open. We say $G$ is \bbb{hereditarily just infinite} if in addition $H$ is just infinite for every open subgroup $H$ of $G$.
\end{defn}

 %Every permutation group can be seen as an abstract group and we will make the convention that any group-theoretical property of a permutation group refers to the group structure, i.e.\! a simple transitive permutation group is a simple group seen as a transitive subgroup of some symmetric group. 
  
In \cite[Theorem~A]{wilson:largehereditarily} the first example of a family of non-(virtually pro-$p$) hereditarily just infinite profinite groups was introduced: the \emph{Wilson groups}. Very few properties of these groups have been discovered, namely: there exist a Wilson group $\mathcal{W}$ such that any finitely generated profinite group can be embedded in $\mathcal{W}$ (see \cite{wilson:largehereditarily}); almost every Wilson group is positively two generated (see \cite{quick:probabilisticgeneration}). %The latter result relies on the classification of finite non-abelian simple groups.

By \cite[Theorem 3]{wilson:justinfinite} a just infinite profinite group is either a \emph{branch group} or it contains an open subgroup isomorphic to the direct product of a finite number of copies of a hereditarily just infinite profinite group. While branch groups have been widely studied in past, not much is known for the second class of groups. The discovery of new properties of Wilson groups will give great insight in the theory of hereditarily just infinite profinite groups. 

In this paper we study a generalisation of Wilson groups that we will call \emph{generalised Wilson groups}. These groups are a special type of \emph{infinitely iterated wreath products} of the groups in $\mS$.

\begin{defn}\label{defn:arbitrary}
 Define the sequence of permutation groups $(\wh{S}_k\le \S(\wh{m}_k))_{k\in \mathbb{N}}$ in the following way: set $\wh{S}_1 = S_1$ and let $\wh{S}_{n+1}$ be the abstract wreath product $S_{n+1} \wr \wh{S}_n$ with a chosen transitive and faithful action on $\wh{m}_k$ points. The groups $\wh{S}_k$ form an inverse system of finite groups and we call their inverse limit $\varprojlim \wh{S}_k$ an \bbb{infinitely iterated wreath product of type $\mS$}. 
\end{defn}

We remark that, for a fixed sequence $\mS$, non-equivalent choices for the action of $\wh{S}_{k}$ in the previous definition lead to non-isomorphic infinitely iterated wreath products.

\begin{defn}
  Let $(k_n)_{n\in \mathbb{N}}$ be an increasing sequence of positive integers. 
  
  A \bbb{generalised Wilson group of type $(\mS,(k_n)_{n\in \mathbb{N}})$} is an infinitely iterated wreath product of type $\mS$, $\varprojlim \wh{S}_k$, where for $n\in \mathbb{N}$ the action of $\wh{S}_{k_n+1} = S_{k_n+ 1}\wr \wh{S}_{k_n}$ is chosen to be the product action of the wreath product. We will write GW group for short. 
%   Define the sequence of permutation groups $\{G_n \le \S(\ol{m}_n)\}_{n\in \mathbb{N}}$ starting from the groups in $\mS$ and the sequence $\{k_n\}_{n\in\mathbb{N}}$ as follows: $\ol{m}_0 = 1$ and $G_0 = \{1\}$; for $n,l\in \mathbb{N}$ and $l\neq k_n-1$, $G_{l}$ is the abstract wreath product $S_{l} \wr G_{l-1}$ with a chosen transitive and faithful action on $\ol{m}_{l}$ points; for $n\in \mathbb{N}$, $G_{k_n-1}$ is the exponentiation $S_{k_n-1} \cwr G_{k_n-2}$. The sequence of finite permutation groups $\{G_{k_n}\}_{n\in \mathbb{N}}$, together with the projections $G_{k_{n+1}} \rightarrow G_{k_n}$, forms an inverse system of finite groups. The inverse limit $\varprojlim G_{k_n}$ is a \bbb{generalised Wilson group of type $(\mS,\{k_n\}_{n\in\mathbb{N}})$}, we will write GW group of type $(\mS,\{k_n\}_{n\in\mathbb{N}})$ for conciseness. 
  % When the choice of the action in the previous construction is not relevant we will denote a GW group obtained from the sequence $\mathcal{S}$ with $\mathrm{GW}(\mathcal{S})$.
  
  An \bbb{infinitely iterated exponentiation of type $\mS$} is the infinitely iterated wreath product of type $\mS$ where every unspecified action is the product action.
\end{defn}

Every infinitely iterated exponentiation of type $\mS$ is a GW group of type $(\mS,(n)_{n\in \mathbb{N}})$. Generalised Wilson groups satisfy the hypotheses of \cite[Theorem~6.2]{reid}, hence they are hereditarily just infinite and not virtually pro-$p$.

% In the previous definition we can take the unspecified actions to be always the exponentiation of the groups in the sequence. 
% \begin{defn}
% Define the sequence of permutation groups $\{\wt{S}_n\le \S(\wt{m}_n)\}_{n\in\mathbb{N}}$ in the following way: $\wt{m}_1 = m_1$, $\wt{S}_1 = S_1$ and $\wt{m}_{n+1} = m_{n+1}^{\wt{m}_n}$, $\wt{S}_{n+1} = S_{n+1}\cwr \wt{S}_n$ for $n\ge 1$. The \bbb{infinitely iterated exponentiation of type $\mS$} is the inverse limit $\varprojlim \wt{S}_n$ of the groups $\wt{S}_{n}$.
% \end{defn}

The goal of this paper is to determine the behaviour of two profinite generation properties in the family of generalised Wilson groups: lower rank and profinite presentability.

In Section \ref{sec:lowerrank} we study the lower rank of an infinitely iterated exponentiation. %A \emph{filtration} of a profinite group $G$ is a descending chain of open subgroups of $G$ which form a base for the topology at the identity.
\begin{defn}
 The \bbb{lower rank} of a profinite group $G$ is the minimal integer $r$ such that $G$ has a base for the neighbourhoods of the identity consisting of $r$-generated subgroups.
\end{defn}

Only very few profinite groups are known to have finite lower rank, amongst them are $p$-adic analytic pro-$p$ groups and $SL_2^1(\mathbb{F}_p[[t]])$, thus it is natural to search for new examples. Our first result will give a new family of profinite groups with finite lower rank. 

For a sequence of integers $(m_k)_{k\in \mathbb{N}}$ we write $\wt{m}_1 = m_1$ and $\wt{m}_{k+1} = m_{k+1}^{\wt{m}_k}$. Let $G$ be a finite group, $\d(G)$ will denote the minimal number of generators of $G$.

\begin{thm}\label{thm:finitelowerrank}
Suppose there is a fixed $r\in \mathbb{N}$ such that $\d(S_{n}^{\wt{m}_{n-1}}) \le r$ for infinitely many $n \in \mathbb{N}$. Then the infinitely iterated exponentiation of type $\mS$ has lower rank at most $r$. 
\end{thm}

% The previous theorem gives a new family of profinite groups with finite lower rank. 

In Section \ref{sec:presentation} we work on the topological finite presentability of infinitely iterated wreath products (see Section~\ref{sec:presentation} for the definition of finite presentability).

As a consequence of \cite[Theorem~A]{quick:probabilisticgeneration}, any infinitely iterated wreath product of type $\mS$ is 2-generated, provided that $m_1>35$. It is then natural to ask whether infinitely iterated wreath products are finitely presentable. The second result of this paper is the following. We denote the Schur multiplier of a finite group $G$ by $M(G)$. 

\begin{thm}\label{thm:nofinpres}
Suppose that the profinite group $\prod_{n\in \mathbb{N}} M(S_n)$ is not topologically finitely generated, then an infinitely iterated wreath product of type $\mS$ is not topologically finitely presentable.
\end{thm}

As a corollary of Theorem \ref{thm:nofinpres} we obtain a sufficient condition for the non-presentability of a generalised Wilson group.

\begin{cor}\label{cor:B}
 Suppose that the profinite group $\prod_{n\in \mathbb{N}} M(S_n)$ is not topologically finitely generated, then a generalised Wilson group of type $(\mS,(k_n)_{n\in \mathbb{N}})$ is not topologically finitely presentable.
\end{cor}

% See Remark \ref{rmk:pres} for a number of cases where the hypothesis of Theorem~\ref{thm:nofinpres} holds. 

\section{Proof of Theorem~\ref{thm:finitelowerrank}}\label{sec:lowerrank}
In this section we discuss the lower rank of infinitely iterated exponentiations. Let $G$ be a finite group, we will write $\d(G)$ for the minimal size of a generating subset of $G$. The following lemma is straightforward.
\begin{lem}\label{lem:minnorm}
 Let $A\le \S(m)$ and $B\le \S(n)$ be two permutation groups. Then $A^n$ is the unique minimal normal subgroup of $A\wr B$ if and only if $A$ is simple and $B$ is transitive.
\end{lem}

If the groups in the sequence $\mS$ ``grow'' sufficiently fast then the infinitely iterated exponentiation obtained from the sequence has finite lower rank. Remember that, for a sequence $(m_k)_{k\in \mathbb{N}}$, we write $\wt{m}_1 = m_1$ and $\wt{m}_{k+1} = m_{k+1}^{\wt{m}_k}$.
% \begin{namedtheorem}[Theorem~A]
%  Let $\mathcal{S} = \{S_k \le \S(m_k)\}_{k \in \mathbb{N}}$ be a sequence of finite non-abelian simple transitive permutation groups. Suppose there is a fixed $r\in \mathbb{N}$ such that $\d(S_{k}^{\wt{m}_{k-1}}) \le r$ for infinitely many $k \in \mathbb{N}$. Then the infinitely iterated exponentiation of type $\mS$ has lower rank at most $r$. 
 \begin{namedproof}[Proof of Theorem \ref{thm:finitelowerrank}]
  Set $G= \varprojlim \wh{S}_k$ where the groups $\wh{S}_k$ are considered all with the product action of the wreath product. Consider the subgroups $N_k= \mathrm{ker}(G \rightarrow \wh{S}_{k})$ for $k \in \mathbb{N}$. It is clear that these subgroups form a base for the topology at the identity and, by Lemma~\ref{lem:minnorm}, they are the only open normal subgroups of $G$. 
  Moreover, $N_{k}/N_{k+1}$ is isomorphic to $S_{k+1}^{\wt{m}_{k}}$ for every $k\ge 1$.
  
  By definition of product action and by Lemma~\ref{lem:minnorm}, $N_k/N_{k+1}$ is the unique minimal normal subgroup of $N_i/N_{k+1}$ every $k\in \mathbb{N}$ and for every $i=1,\ldots, k$. 
%   , the action of $S_k^{\wt{m}_k}$ as a subgroup of $\wh{S}_k$ is transitive. Using Lemma~\ref{lem:minnorm}, it is now easy to show that $N_k/N_{k+1}$ is the unique minimal normal subgroup of $N_i/N_{k+1}$ every $k\in \mathbb{N}$ and for every $i=1,\ldots, k$. 
Repeated applications of \cite[Theorem~1.1]{MR1492976} yield $\d(N_{k-1}) = \mathrm{d}(N_{k-1}/N_k)$. By hypothesis, $\mathrm{d}(S_{k}^{\wt{m}_{k-1}}) \le r$ for infinitely many $k$ and $\{N_k \vert \d(N_k)\le r\}$ is the required base of $G$. \qed
 \end{namedproof}
% \end{namedtheorem}
\begin{rmk}
 The hypotheses of Theorem~A are satisfied, with $r=2$, by the sequences $\mS = (\mathrm{PSL}_2(p_n)\le \S(p_n+1))_{n\in \mathbb{N}}$ such that $\mathrm{PSL}_2(p_n)$ acts on the projective line over $\mathbb{F}_{p_n}$ and $(p_n)_{n\in \mathbb{N}}$ is any sequence of primes satisfying
\[
  p_n \ge \frac{1}{4}(p_{n-1}+1)(p_{n-1}^2 - 2p_{n-1}-1)-2.
\]
This follows from the calculation of the \emph{Eulerian function} for $\mathrm{PSL}_2(p)$ (see \cite{hall:eulerianfunctions}). In particular, infinitely iterated exponentiations of these sequences have lower rank 2. 
\end{rmk}
%  We believe that the behaviour of the lower rank in GW groups is various. 
 
 We conjecture that there exists a GW group with infinite lower rank. One strong candidate for this is the infinitely iterated exponentiation of a constant sequence, i.e.\! $S_{k} = S$ for $k\in \mathbb{N}$. We are also convinced that there are GW groups of arbitrary finite lower rank, but a proof of this result has to involve an accurate study of the subgroup structure of GW groups.
% \begin{conj}
%  There exists a generalised Wilson group of infinite lower rank.
% \end{conj}
\begin{conj}
 There exists a generalised Wilson group of lower rank $r$ for every $r \in \mathbb{N}\cup \{\infty\}$.
\end{conj}
A positive answer to this conjecture would produce interesting examples as the only known family of profinite groups of arbitrary finite lower rank are $p$-adic analytic pro-$p$ groups, where the lower rank coincides with the number of generators of the associated $p$-adic Lie algebra (see \cite{powerful}).

\section{Proof of Theorem~\ref{thm:nofinpres}}\label{sec:presentation}
Let $H$ be a finite perfect group and fix a surjective homomorphism $F\rightarrow H$ with kernel $R$ from an appropriate free group $F$. The \emph{Schur multiplier} of $H$ is the finite group $R/[F,R]$. Let $G$ be a profinite group, we will write $\d(G)$ for the minimal size of a subset of $G$ which generates a dense subgroup of $G$.
% \begin{defn}
%  Let $G$ be a profinite group and 
Let $N$ be a closed normal subgroup of $G$ and let $\mathcal{R}$ be a subset of $N$, we say that $N$ is \emph{topologically normally generated in $G$ by $\mathcal{R}$} if the $G$-conjugates of the elements of $\mathcal{R}$ generate a dense subgroup of $N$.
% \[
%   N =\ol{ \langle g^{-1} s g \vvert g\in G, s\in \mathcal{R} \rangle}.
% \]
% \end{defn}
 
Let $G$ be a topologically $d$-generated profinite group. We can define a continuous epimorphism $F\rightarrow G$ from the free profinite group $F= \widehat{F}_{d}$; the kernel of this epimorphism is a closed normal subgroup $R$ of $F$. Let $\mathcal{S}$ be a set of topological generators for $F$ and let $R$ be topologically normally generated by a subset $\mathcal{R}$ of $R$, then these give us a \bbb{profinite presentation} of $G$ and we write $G = \langle \mathcal{S} \vert \mathcal{R} \rangle$. We say that a finitely generated profinite group $G$ is \bbb{topologically finitely presentable} if there exists a presentation $G = \langle \mathcal{S} \vert \mathcal{R} \rangle$ of $G$ such that $\mathcal{S}$ and $\mathcal{R}$ are finite.

 The following lemma is an application of \cite[Theorem~3]{MR0414723}.
\begin{lem}\label{lem:read}
 Suppose that $\prod_{n\in \mathbb{N}} M(S_n)$ is not topologically finitely generated.  Then the sequence $(\d(M(\wh{S}_n)))_{n\in \mathbb{N}}$ is unbounded. 
\begin{proof}
  By \cite[Theorem~3]{MR0414723}, if $A$ is a perfect group and $B$ is a perfect permutation group, then $M(A\wr B) \cong M(A)\times M(B)$. Therefore $M(\wh{S}_n) = \prod_{k=1}^n M(S_k)$, the claim follows.
\end{proof}
\end{lem}

We now apply Lemma \ref{lem:read} to the proof of Theorem~\ref{thm:nofinpres}. %give sufficient conditions for an infinitely iterated wreath product to be not finitely presentable.

% \begin{namedtheorem}[Theorem~\ref{thm:nofinpres}]
% Suppose that the profinite group $\prod_{n\in \mathbb{N}} M(S_n)$ is not topologically finitely generated, then an infinitely iterated wreath product of type $\mS$ is not topologically finitely presentable.
\begin{namedproof}[Proof of Theorem~\ref{thm:nofinpres}]
 %Let $G$ be a GW group of type $\mS$. If $G$ is not topologically finitely generated, then $G$ is not finitely presentable.
 Let $G$ be a $d$-generated infinitely iterated wreath product of type $\mS$. Let $F=\widehat{F}_d$ be the profinite free group of rank $d$, $\varphi: F\rightarrow G$ a continuous epimorphism, $R = \mathrm{ker}\varphi$ and let $N$ be an open normal subgroup of $F$. Now, $NR$ is an open normal subgroup of $F$ that contains $R$, thus $F/NR$ is isomorphic to a continuous quotient of $G$. By Lemma \ref{lem:minnorm}, the only open normal subgroups of $G$ are kernels of the inverse limit projections from $G$ to $G_n$ for some integer $n$. Therefore, any continuous quotient of $G$ is isomorphic, as an abstract group, to an iterated wreath product.
 
 The number of relations of $G$ in the chosen presentation is the (possibly infinite) number $r(G)$ of normal generators for $R$ as a subgroup of $F$. The quotient $R/[R,R]$ of $R$ is abelian and $R/[F,R]$ is a quotient of the latter, hence $r(G) \ge \mathrm{d}(R/[R,R]) \ge \mathrm{d}(R/[F,R])$.
 %, since both $N$ and $R$ are normal in $F$ we have $M= N[F,R]$. 

Set $M= [F,NR]$. $R/(M\cap R)$ is a quotient of $R/[F,R]$, hence $\d(R/[F,R]) \ge \d(R/(M\cap R))$. $R/(M\cap R)$ is isomorphic to $NR/M$ and $NR/M = NR/[F,NR]$ is the Schur multiplier of an iterated wreath product, as shown in the previous paragraph. Therefore $r(G)\ge \d(NR/[F,NR])$. By Lemma \ref{lem:read} and by hypothesis, the last quantity is unbounded as $N$ ranges between all open normal subgroups of $F$. Thus $G$ cannot be finitely presentable. \qed
\end{namedproof}
% \end{namedtheorem}

\begin{rmk}\label{rmk:pres}
A sequence of finite non-abelian simple groups $(S_k)_{k\in \mathbb{N}}$ such that a fixed prime $p$ divides $M(S_n)$ for infinitely many $n$ satisfies the hypotheses of Theorem~\ref{thm:nofinpres}. Let $\mathcal{C}$ be the constant sequence $(\mathrm{Alt}(36)\le \S(36))_{k\in \mathbb{N}}$, then every infinitely iterated wreath product of type $\mathcal{C}$ is finitely generated by \cite{quick:probabilisticgeneration}, but it is not finitely presentable by Theorem~\ref{thm:nofinpres}. In fact $M(\mathrm{Alt}(36))$ has order two.
% In general, if $S$ is a finite non-abelian simple group with non-trivial Schur multiplier and $S$ is not a Chevalley group of type $A_n$ or a Steinberg group of type $\phantom{}^2A_n$, then $\norm{M(S)}$ is divisible by 2 or 3 (possibly both).
\end{rmk}
We believe that no generalised Wilson group is finitely presentable, but we do not have a proof of this general statement.

\begin{conj}
Every generalised Wilson group is not finitely presentable.
\end{conj}

\section*{Acknowledgements}
This work was carried out as part of my Ph.D studies at Royal Holloway University of London under the supervision of Dr.~Yiftach Barnea, I would like to thank him for introducing me to the problem and his continuous help and support. I also wish to thank Eugenio Giannelli and Benjamin Klopsch for the most useful discussions, in particular I am grateful to the latter for the suggestion of Theorem~\ref{thm:nofinpres}.

\bibliographystyle{plain}

\begin{thebibliography}{10}

\bibitem{wilson:justinfinite}
 R.~Grigorchuk.
\newblock Just infinite branch groups.
\newblock {\em New horizons in pro-{$p$} groups}, {\em Progr. Math.} 184:121--179.
\newblock Birkh\"auser Boston, Boston, MA, 2000.

\bibitem{hall:eulerianfunctions}
P.~Hall.
\newblock Eulerian functions of a group.
\newblock {\em Quart. J. Math.}, 7:134--151, 1936.

\bibitem{powerful}
A.~Lubotzky and A.~Mann.
\newblock Powerful {$p$}-groups. {II}. {$p$}-adic analytic groups.
\newblock {\em J. Algebra}, 105(2):506--515, 1987.

\bibitem{MR1492976}
A.~Lucchini and F.~Menegazzo.
\newblock Generators for finite groups with a unique minimal normal subgroup.
\newblock {\em Rend. Sem. Mat. Univ. Padova}, 98:173--191, 1997.

\bibitem{quick:probabilisticgeneration}
M.~Quick.
\newblock Probabilistic generation of wreath products of non-abelian finite
  simple groups. {II}.
\newblock {\em Internat. J. Algebra Comput.}, 16(3):493--503, 2006.

\bibitem{MR0414723}
E.~W.~Read.
\newblock On the {S}chur multiplier of a wreath product.
\newblock {\em Illinois J. Math.}, 20(3):456--466, 1976.

 \bibitem{reid}
C.~D.~Reid.
 \newblock Inverse system characterizations of the (hereditarily) just infinite property in profinite groups.
 \newblock {\em Bull. Lond. Math. Soc.}, 44(3):413--425, 2012. 

\bibitem{wilson:largehereditarily}
J.~S.~Wilson.
\newblock Large hereditarily just infinite groups.
\newblock {\em J. Algebra}, 324(2):248--255, 2010.

\end{thebibliography}
\def\cprime{$'$}

\section*{Author information}
Matteo Vannacci\\
Royal Holloway, University of London\\
Egham, Surrey, TW200EX\\
United Kingdom.\\
E-mail: \texttt{vannacci.m@gmail.com}

\end{document}